\documentclass[a4paper,11pt]{article}

\usepackage{amsfonts,amssymb,amsmath,amsthm,latexsym,epsf}

\usepackage[english]{babel}
\makeatletter
\def\@captype{figure}
\makeatother

\makeatletter
\@addtoreset{equation}{section}
\makeatother

\newcounter{enunciato}[section]

\newtheorem{ittheorem}{Theorem}
\newtheorem{itlemma}{Lemma}
\newtheorem{itconjecture}{Conjecture}

\newenvironment{theorem}{\addtocounter{enunciato}{1}
\begin{ittheorem}}{\end{ittheorem}}

\newenvironment{conjecture}{\addtocounter{enunciato}{1}
\begin{itconjecture}}{\end{itconjecture}}

\begin{document}
\title{Intermittency on catalysts}

\author{\renewcommand{\thefootnote}{\arabic{footnote}}
J.\ G\"artner
\footnotemark[1]
\\
\renewcommand{\thefootnote}{2, 3}
F.\ den Hollander
\footnotemark
\\
\renewcommand{\thefootnote}{\arabic{footnote}}
G.\ Maillard
\footnotemark[4]
}

\footnotetext[1]{
Institut f\"ur Mathematik, Technische Universit\"at Berlin,
Strasse des 17.\ Juni 136, D-10623 Berlin, Germany,
{\sl jg@math.tu-berlin.de}
}
\footnotetext[2]{
Mathematical Institute, Leiden University, P.O.\ Box 9512,
2300 RA Leiden, The Netherlands,
{\sl denholla@math.leidenuniv.nl}
}
\footnotetext[3]{
EURANDOM, P.O.\ Box 513, 5600 MB Eindhoven, The Netherlands
}
\footnotetext[4]{
Institut de Math\'ematiques, \'Ecole Polytechnique F\'ed\'erale
de Lausanne, CH-1015 Lausanne, Switzerland,
{\sl gregory.maillard@epfl.ch}
}

\maketitle

\begin{abstract}

The present paper provides an overview of results obtained in four
recent papers by the authors. These papers address the problem of
intermittency for the Parabolic Anderson Model in a \emph{time-dependent
random medium}, describing the evolution of a ``reactant'' in the
presence of a ``catalyst''. Three examples of catalysts are considered:
(1) independent simple random walks; (2) symmetric exclusion process;
(3) symmetric voter model. The focus is on the annealed Lyapunov exponents,
i.e., the exponential growth rates of the successive moments of the
reactant. It turns out that these exponents exhibit an interesting
dependence on the dimension and on the diffusion constant.

\vspace{0.5cm}\noindent

\vspace{0.5cm}\noindent
* Invited paper to appear in a Festschrift in honour of Heinrich von
Weizs\"acker, on the occasion of his 60th birthday, to be published
by Cambridge University Press.
\end{abstract}

\newpage


\section{The Parabolic Anderson Model}
\label{S1}

\subsection{Motivation}
\label{S1.1}

The {\em Parabolic Anderson Model\/} is the partial differential equation
\begin{equation}
\label{pA}
\frac{\partial}{\partial t}u(x,t) = \kappa\Delta u(x,t) + \gamma \xi(x,t)u(x,t),
\qquad x\in\mathbb{Z}^d,\,t\geq 0,
\end{equation}
for the $\mathbb{R}$-valued random field
\begin{equation}
\label{rfu}
u = \{u(x,t) \colon\,x\in\mathbb{Z}^d,\,t \geq 0\},
\end{equation}
where $\kappa\in [0,\infty)$ is the diffusion constant, $\gamma\in [0,\infty)$ is
the coupling constant, $\Delta$ is the discrete Laplacian, acting on $u$ as
\begin{equation}
\label{dL}
\Delta u(x,t) = \sum_{{y\in\mathbb{Z}^d} \atop {\|y-x\|=1}} [u(y,t)-u(x,t)]
\end{equation}
($\|\cdot\|$ is the Euclidian norm), while
\begin{equation}
\label{rf}
\xi = \{\xi(x,t) \colon\,x\in\mathbb{Z}^d,\,t \geq 0\}
\end{equation}
is an $\mathbb{R}$-valued random field that evolves with time and that drives
the equation. As initial condition for (\ref{pA}) we take
\begin{equation}
\label{ic}
u(\cdot\,,0) \equiv 1.
\end{equation}

One interpretation of (\ref{pA}) and (\ref{ic}) comes from {\em population dynamics}.
Consider a spatially homogeneous system of two types of particles, $A$ (catalyst)
and $B$ (reactant), subject to:
\begin{itemize}
\item[(i)]
$A$-particles evolve autonomously, according to a prescribed stationary dynamics
given by the $\xi$-field, with $\xi(x,t)$ denoting the number of $A$-particles
at site $x$ at time $t$;
\item[(ii)]
$B$-particles perform independent simple random walks with jump rate $2d\kappa$
and split into two at a rate that is equal to $\gamma$ times the number of 
$A$-particles present at the same location;
\item[(iii)]
the initial density of $B$-particles is $1$.
\end{itemize}
Then
\begin{equation}
\label{uint}
\begin{array}{lll}
u(x,t) &=& \hbox{the average number of $B$-particles at site $x$ at time $t$}\\
       && \hbox{conditioned on the evolution of the $A$-particles}.
\end{array}
\end{equation}
It is possible to add that $B$-particles die at rate $\delta\in (0,\infty)$. This
amounts to the trivial transformation
\begin{equation}
\label{udel}
u(x,t) \to u(x,t)e^{-\delta t}.
\end{equation}

What makes (\ref{pA}) particularly interesting is that the two terms in the right-hand 
side \emph{compete with each other\/}: the diffusion (of $B$-particles) described by 
$\kappa\Delta$ tends to make $u$ flat, while the branching (of $B$-particles caused 
by $A$-particles) described by $\xi$ tends to make $u$ irregular.

\subsection{Intermittency}
\label{S1.2}

We will be interested in the presence or absence of {\em intermittency}.
Intermittency means that for large $t$ the branching dominates, i.e.,
the $u$-field develops sparse high peaks in such a way that $u$ and its moments
are each dominated by their own collection of peaks (see G\"artner and K\"onig
\cite{garkon05}, Section 1.3). In the {\em quenched\/} situation, i.e.,
conditional on $\xi$, this geometric picture of intermittency is well understood
for several classes of \emph{time-independent\/} random potentials $\xi$ (see
e.g.\ Sznitman \cite{sz98} for Poisson clouds and G\"artner, K\"onig and Molchanov
\cite{garkonmol05} for i.i.d.\ potentials with double-exponential and heavier
upper tails; G\"artner and K\"onig \cite{garkon05} provides an overview). For 
\emph{time-dependent\/} random potentials $\xi$, however, such a geometric
picture is not yet available. Instead one restricts attention to understanding
the phenomenon of intermittency indirectly by comparing the successive
{\em annealed\/} Lyapunov exponents
\begin{equation}
\label{lyapdef}
\lambda_p = \lim_{t\to\infty} \Lambda_p(t), \qquad p\in\mathbb{N},
\end{equation}
with
\begin{equation}
\label{lyappredef}
\Lambda_p(t) = \frac{1}{t}\log\mathbb{E}\left([u(0,t)]^p\right)^{1/p},
\qquad p \in \mathbb{N},\,t > 0,
\end{equation}
where $\mathbb{E}$ denotes expectation w.r.t.\ $\xi$. One says that the solution $u$
is {\em $p$-intermittent\/} if
\begin{equation}
\label{lambstrict}
\lambda_p > \lambda_{p-1},
\end{equation}
and {\em intermittent\/} if (\ref{lambstrict}) holds for all $p\in\mathbb{N}
\setminus\{1\}$.

Carmona and Molchanov~\cite{carmol94} succeeded to investigate the annealed Lyapunov
exponents, and to obtain the qualitative picture of intermittency (in terms of these
exponents), for potentials of the
form
\begin{equation}
\label{whitenoisepot}
\xi(x,t) = \dot{W}_x(t),
\end{equation}
where $\{W_x(t)\colon\,x\in\mathbb{Z}^d,\,t\geq 0\}$ denotes a collection of
independent Brownian motions. (In this case, (\ref{pA}) corresponds to an infinite
system of coupled It\^o-diffusions.) They showed that for $d=1,2$ intermittency
holds for all $\kappa$, whereas for $d\geq 3$ $p$-intermittency holds if and
only if the diffusion constant $\kappa$ is smaller than a critical threshold
$\kappa_p=\kappa_p(d,\gamma)$ tending to infinity as $p\to\infty$. They also 
studied the asymptotics of the quenched Lyapunov exponent in the limit as 
$\kappa\downarrow 0$, which turns out to be singular. Subsequently, the latter 
was more thoroughly investigated in papers by Carmona, Molchanov and 
Viens~\cite{carmolvi96}, Carmona, Koralov and Molchanov~\cite{carkormol01}, 
and Cranston, Mountford and Shiga~\cite{cramoshi02}.

In Sections~\ref{S2}--\ref{S4} we consider three different choices for $\xi$,
namely:
\begin{itemize}
\item[(1)] Independent Simple Random Walks.
\item[(2)] Symmetric Exclusion Process.
\item[(3)] Symmetric Voter Model.
\end{itemize}
For each of these examples we study the annealed Lyapunov exponents as a function
of $d$, $\kappa$ and $\gamma$. Because of their \emph{non-Gaussian\/} and 
\emph{non-independent\/} spatial structure, these examples require techniques 
different from those developed for (\ref{whitenoisepot}). Example (1) was studied 
earlier in Kesten and Sidoravicius~\cite{kessid03}. We describe their work in 
Section~\ref{S2.2}.

By the Feynman-Kac formula, the solution of (\ref{pA}) and (\ref{ic}) reads
\begin{equation}
\label{fey-kac1}
u(x,t) = E_{\,x}\left(\exp\left[\gamma \int_0^t ds\,\,
\xi\left(X^\kappa(s),t-s\right)\right]\right),
\end{equation}
where $X^\kappa$ is simple random walk on $\mathbb{Z}^d$ with step rate $2d\kappa$
and $E_{\,x}$ denotes expectation with respect to $X^\kappa$ given $X^\kappa(0)=x$.
This formula shows that understanding intermittency amounts to studying the
{\em large deviation behavior\/} of a random walk {\em sampling\/} a time-dependent
random field.


\section{Independent Simple Random Walks}
\label{S2}

In this section we consider the case where $\xi$ is a Poisson field of
{\em Independent Simple Random Walks\/} (ISRW). We first describe the
results obtained in Kesten and Sidoravicius~\cite{kessid03}. After that
we describe the refinements of these results obtained in G\"artner and
den Hollander~\cite{garhol06}.

\subsection{Model}
\label{S2.1}

ISRW is the Markov process with state space
\begin{equation}
\label{statIRW}
\Omega = (\mathbb{N} \cup \{0\})^{\mathbb{Z}^d}
\end{equation}
whose generator acts on cylindrical functions $f$ as
\begin{equation}
\label{genIRW}
(Lf)(\eta) = \frac{1}{2d}\,\sum_{(x,y)} \eta(x)
[f(\eta^{x \curvearrowright y})-f(\eta)],
\end{equation}
where the sum runs over oriented bonds between neighboring sites, and
\begin{equation}
\label{irwpro}
\eta^{x \curvearrowright y}(z) =
\begin{cases}
\eta(z)    &\text{if } z \neq x,y,\\
\eta(x)-1  & \text{if } z=x,\\
\eta(y)+1  & \text{if } z=y,
\end{cases}
\end{equation}
i.e., $\eta^{x \curvearrowright y}$ is the configuration obtained from $\eta$
by moving a particle from $x$ to $y$. We choose $\xi(\cdot\,,0)$ according to
the Poisson product measure with density $\rho \in (0,\infty)$, i.e., initially
each site carries a number of particles that is Poisson distributed with mean
$\rho$. For this choice, the $\xi$-field is stationary and reversible in time 
(see Kipnis and Landim~\cite{kiplan99}).

Under ISRW, particles move around independently as simple random walks, stepping
at rate $1$ and choosing from neighboring sites with probability $1/2d$ each.

\subsection{Main theorems}
\label{S2.2}

Kesten and Sidoravicius~\cite{kessid03} proved the following. They considered the
language of $A$-particles and $B$-particles from population dynamics, as mentioned
in Section \ref{S1.1}, and included a death rate $\delta \in [0,\infty)$ for the
$B$-particles (recall (\ref{udel})).

\begin{itemize}
\item[(1)]
If $d=1,2$, then -- for any choice of the parameters -- the average number of
$B$-particles per site tends to infinity at a rate that is faster than
exponential.
\item[(2)]
If $d \geq 3$, then -- for $\gamma$ sufficiently small and $\delta$ sufficiently
large -- the average number of $B$-particles per site tends to zero exponentially 
fast.
\item[(3)]
If $d \geq 1$, then -- conditional on the evolution of the $A$-particles -- there 
is a phase transition: for small $\delta$ the $B$-particles locally survive, while 
for large $\delta$ they become locally extinct.
\end{itemize}
Properties (1) and (2) -- which are annealed results -- are implied by
Theorems~\ref{LyadiIRW} and \ref{LyapropIRW} below, while property (3)
-- which is a quenched result -- is not. The main focus of \cite{kessid03}
is on survival versus extinction. The approach in \cite{kessid03}, being
based on path estimates rather than on the Feynman-Kac representation,
produces cruder results, but it is more robust against variations of the
dynamics.

In G\"artner and den Hollander~\cite{garhol06} the focus is on the annealed
Lyapunov exponents. Theorems~\ref{LyaexistIRW}--\ref{LyapropIRW} below are
taken from that paper.

\begin{theorem}
\label{LyaexistIRW}
Let $d\geq 1$, $\rho,\gamma\in (0,\infty)$ and $p\in\mathbb{N}$.\\
(i) For all $\kappa\in [0,\infty)$, the limit in {\rm (\ref{lyapdef})} exist.\\
(ii) If $\lambda_p(0)<\infty$, then $\kappa\to\lambda_p(\kappa)$
is finite, continuous, non-increasing and convex on $[0,\infty)$.
\end{theorem}

Let $p_t(x,y)$ denote the probability that simple random walk stepping at
rate 1 moves from $x$ to $y$ in time $t$. Let
\begin{equation}
\label{Gddef}
G_d = \int_0^\infty p_t(0,0)\,dt
\end{equation}
be the Green function at the origin of simple random walk.

\begin{theorem}
\label{LyadiIRW}
Let $d\geq 1$, $\rho,\gamma\in (0,\infty)$ and $p\in\mathbb{N}$. Then, for
all $\kappa\in [0,\infty)$, $\lambda_p(\kappa)<\infty$ if and only if
$p<1/G_d\gamma$.
\end{theorem}

\noindent
It can be shown that if $p > 1/G_d\gamma$, then $\Lambda_p(t)$ in (\ref{lyappredef})
grows exponentially fast with $t$, i.e., the $p$-th moment of $u(0,t)$ grows
double exponentially fast with $t$. The constant in the exponent can be computed.

In the regime $p<1/G_d\gamma$, $\kappa\mapsto\lambda_p(\kappa)$ has the following
behavior (see Fig.~\ref{fig-IRWasymp}):

\begin{theorem}
\label{LyapropIRW}
Let $d\geq 1$, $\rho,\gamma\in (0,\infty)$ and $p\in\mathbb{N}$ such that 
$p < 1/G_d\gamma$.\\
(i) $\kappa\mapsto\lambda_p(\kappa)$ is continuous, strictly decreasing and
convex on $[0,\infty)$.\\
(ii) For $\kappa=0$,
\begin{equation}
\label{lam0}
\lambda_p(0) = \rho\gamma\,\,\frac{(1/G_d)}{(1/G_d)-p\gamma}.
\end{equation}
(iii) For $\kappa\to\infty$,
\begin{equation}
\label{limlamb*IRW}
\lim_{\kappa\to\infty} 2d\kappa[\lambda_p(\kappa)-\rho\gamma]
= \rho\gamma^2 G_d + 1_{d=3}\,(2d)^3(\rho\gamma^2 p)^2\,{\cal P}_3
\end{equation}
with
\begin{equation}
\label{poldef}
{\cal P}_3 = \sup_{ {f \in H^1(\mathbb R^3)} \atop {\|f\|_2=1} }
\Big[\, \int_{\mathbb{R}^3} dx\,|f(x)|^2 \int_{\mathbb{R}^3} dy\,|f(y)|^2
\,\frac{1}{4\pi\|x-y\|}
- \int_{\mathbb{R}^3} dx\,|\nabla f(x)|^2\,\Big].
\end{equation}
\end{theorem}

\vskip 0.8truecm

\begin{center}
\setlength{\unitlength}{0.3cm}
\begin{picture}(20,10)(12,0)
\put(24,-4){\line(18,0){18}}
\put(24,-4){\line(0,14){14}}
{\thicklines
\qbezier(24,8)(26,6.6)(28,5.5)
\qbezier(24,6)(26,4.8)(28,3.9)
\qbezier(24,4)(26,3)(28,2.4)
}
\qbezier[25](36,0.8)(38,0.5)(41,0.4)
\qbezier[60](24,0)(33,0)(42,0)
\put(22,-.5){$\rho\gamma$}
\put(23,-4.5){$0$}
\put(24,8){\circle*{.35}}
\put(24,6){\circle*{.35}}
\put(24,4){\circle*{.35}}
\put(31,2.5){{\bf ?}}
\put(43,-4.1){$\kappa$}
\put(22.5,10.7){$\lambda_p(\kappa)$}
\put(30,8){$d \geq 4$}
\put(0,-4){\line(18,0){18}}
\put(0,-4){\line(0,14){14}}
{\thicklines
\qbezier(0,8)(2,6.6)(4,5.3)
\qbezier(0,6)(2,4.8)(4,3.8)
\qbezier(0,4)(2,3)(4,2.25)
}
{\thicklines
\qbezier(11,2.0)(13,1.4)(17,0.85)
\qbezier(11,1.4)(13,0.9)(17,0.55)
\qbezier(11,0.7)(13,0.4)(17,0.25)
}
\qbezier[80](0,0)(9,0)(18,0)
\put(-2,-.5){$\rho\gamma$}
\put(-1,-4.5){$0$}
\put(0,8){\circle*{.35}}
\put(0,6){\circle*{.35}}
\put(0,4){\circle*{.35}}
\put(7,2.5){{\bf ?}}
\put(19,-4.1){$\kappa$}
\put(-1.5,10.7){$\lambda_p(\kappa)$}
\put(6,8){$d=3$}
\end{picture}
\vspace{1.5cm}
\end{center}
\caption{$\kappa\mapsto\lambda_p(\kappa)$ for $p=1,2,3$ when $p<1/G_d\gamma$ for simple
random walk in $d=3$ and $d\geq 4$.}\label{fig-IRWasymp}

\vskip 0.8truecm

\subsection{Discussion}
\label{S2.3}

Theorem \ref{LyadiIRW} says that if the catalyst is driven by a recurrent random
walk ($G_d=\infty$), then it can pile up near the origin and make the reactant
grow at an unbounded rate, while if the catalyst is driven by a transient random
walk ($G_d<\infty$), then small enough moments of the reactant grow at a finite
rate. We refer to this dichotomy as the {\em strongly catalytic}, respectively,
the {\em weakly catalytic\/} regime.

Theorem \ref{LyapropIRW}(i) shows that, even in the weakly catalytic regime, some
degree of {\em clumping\/} of the catalyst occurs, in that the growth rate of the 
reactant is $>\rho\gamma$, the average medium growth rate. As the diffusion constant 
$\kappa$ of the reactant increases, the effect of the clumping of the catalyst on 
the reactant gradually diminishes, and the growth rate of the reactant gradually 
decreases to $\rho\gamma$.

Theorem \ref{LyapropIRW}(ii) shows that, again in the weakly catalytic regime, if
the reactant stands still, then the system is intermittent. Apparently, the
successive moments of the reactant are sensitive to \emph{successive degrees
of clumping}. By continuity, intermittency persists for small $\kappa$.

Theorem \ref{LyapropIRW}(iii) shows that all Lyapunov exponents decay to 
$\rho\gamma$ as $\kappa\to\infty$ in the same manner when $d \geq 4$ but 
not when $d=3$. In fact, in $d=3$ intermittency persists for large $\kappa$. 
It remains open whether the same is true for $d\geq 4$. To decide the latter, 
we need a finer asymptotics for $d\geq 4$. A large diffusion constant of the 
reactant hampers localization of the reactant around regions where the catalyst 
clumps, but it is not a priori clear whether this is able to destroy intermittency 
for $d\geq 4$. We conjecture:

\begin{conjecture}
\label{cj1IRW}
In $d=3$, the system is intermittent for all $\kappa\in [0,\infty)$.
\end{conjecture}

\begin{conjecture}
\label{cj2IRW}
In $d\geq 4$, there exists a strictly increasing sequence $0<\kappa_2<\kappa_3
<\ldots$ such that for $p=2,3,\ldots$ the system is $p$-intermittent if and only
if $\kappa \in [0,\kappa_p)$.
\end{conjecture}

\noindent
In words, we conjecture that in $d=3$ the curves in Fig.~\ref{fig-IRWasymp} never
merge, whereas for $d\geq 4$ the curves merge successively.

What is remarkable about the scaling of $\lambda_p(\kappa)$ as $\kappa\to\infty$
in (\ref{limlamb*IRW}) is that $\mathcal{P}_3$ is the variational problem for
the so-called {\em polaron model\/}. Here, one considers the quantity
\begin{equation}
\label{polintdef}
\theta(t;\alpha) = \frac{1}{\alpha^2t} \log E_0\left(\exp\left[\alpha
\int_0^t ds \int_s^t du\,\,\frac{e^{-(u-s)}}{|\beta(u)-\beta(s)|}\right]\right),
\end{equation}
where $\alpha>0$ and $(\beta(t))_{t\geq 0}$ is standard Brownian motion on
$\mathbb{R}^3$ starting at $\beta(0)=0$. Donsker and Varadhan~\cite{dova83}
proved that
\begin{equation}
\label{polaronlim}
\lim_{\alpha\to\infty} \lim_{t\to\infty} \theta(t;\alpha)
= 4\sqrt{\pi}\,\mathcal P_3.
\end{equation}
Lieb~\cite{li77} proved that (\ref{poldef}) has a unique maximizer modulo
translations and that the centered maximizer is radially symmetric,
radially non-increasing, strictly positive and smooth. A deeper analysis
shows that the link between the scaling of $\lambda_p(\kappa)$ for $\kappa\to\infty$
and the scaling of the polaron for $\alpha\to\infty$ comes from {\em moderate}
deviation behavior of $\xi$ and {\em large} deviation behavior of the occupation
time measure of $X^\kappa$ in (\ref{fey-kac1}). For details we refer to G\"artner 
and den Hollander~\cite{garhol06}.


\section{Symmetric Exclusion Process}
\label{S3}

In this section we consider the case where $\xi$ is the \emph{Symmetric Exclusion
Process} (SEP) in equilibrium. We summarize the results obtained in G\"artner, 
den Hollander and Maillard \cite{garholmai07}, \cite{garholmai07pra}.

\subsection{Model}
\label{S3.1}

Let $p\colon\mathbb{Z}^d\times\mathbb{Z}^d\to [0,1]$ be the transition kernel
of an irreducible symmetric random walk. SEP is the Markov process with state
space
\begin{equation}
\label{statSE}
\Omega = \{0,1\}^{\mathbb{Z}^d}
\end{equation}
whose generator $L$ acts on cylindrical functions $f$ as
\begin{equation}
\label{expro1}
(Lf)(\eta) = \sum_{\{x,y\}\subset\mathbb{Z}^d} p(x,y)\,
\left[f\left(\eta^{x,y}\right)-f(\eta)\right],
\end{equation}
where the sum runs over unoriented bonds between any pair of sites, and
\begin{equation}
\label{expro2}
\eta^{x,y}(z) =
\begin{cases}
\eta(z)  &\text{if } z \neq x,y,\\
\eta(y)  & \text{if } z=x,\\
\eta(x)  & \text{if } z=y.
\end{cases}
\end{equation}
In words, the states of $x$ and $y$ are interchanged along the bond $\{x,y\}$
at rate $p(x,y)$. We choose $\xi(\cdot\,,0)$ according to the Bernoulli product
measure with density $\rho\in (0,1)$. For this choice, the $\xi$-field is
stationary and reversible in time (see Liggett~\cite{lig85}).

Under SEP, particles move around independently according to the symmetric
random walk transition kernel $p(\cdot,\cdot)$, but subject to the restriction
that no two particles can occupy the same site. A special case is simple
random walk
\begin{equation}
\label{SRW}
p(x,y) =
\begin{cases}
\frac{1}{2d} & \text{if } \|x-y\|=1,\\
0            & \text{otherwise}.
\end{cases}
\end{equation}

\subsection{Main theorems}
\label{S3.2}

\begin{theorem}
\label{LyaexistSE}
Let $d\geq 1$, $\rho\in (0,1)$, $\gamma\in (0,\infty)$ and $p\in\mathbb{N}$.\\
(i) For all $\kappa\in [0,\infty)$, the limit in {\rm (\ref{lyapdef})} exists
and is finite.\\
(ii) On $[0,\infty)$, $\kappa\to\lambda_p(\kappa)$ is continuous, non-increasing
and convex.
\end{theorem}

The following dichotomy holds (see Fig.~\ref{fig-SEqual}):

\begin{theorem}
\label{LyalowSE}
Let $d\geq 1$, $\rho\in (0,1)$, $\gamma\in (0,\infty)$ and $p\in\mathbb{N}$.\\
(i) If $p(\cdot,\cdot)$ is recurrent, then $\lambda_p(\kappa)=\gamma$ for all $\kappa
\in [0,\infty)$.\\
(ii) If $p(\cdot,\cdot)$ is transient, then $\rho\gamma<\lambda_p(\kappa)<\gamma$ for 
all $\kappa\in [0,\infty)$. Moreover, $\kappa\mapsto\lambda_p(\kappa)$ is strictly
decreasing with $\lim_{\kappa\to\infty}\lambda_p(\kappa)=\rho\gamma$. Furthermore,
$p\mapsto\lambda_p(0)$ is strictly increasing.
\end{theorem}

For transient simple random walk, $\kappa\mapsto\lambda_p(\kappa)$ has
the following behavior (similar as in Fig.~\ref{fig-IRWasymp}):

\begin{theorem}
\label{LyahighlimSE}
Let $d \geq 3$, $\rho\in (0,1)$, $\gamma\in (0,\infty)$ and $p\in\mathbb{N}$. Assume 
{\rm (\ref{SRW})}. Then
\begin{equation}
\label{limlamb*SE}
\lim_{\kappa\to\infty} 2d\kappa[\lambda_p(\kappa)-\rho\gamma]
= \rho(1-\rho)\gamma^2G_d + 1_{\{d=3\}}\,(2d)^3[\rho(1-\rho)\gamma^2p]^2\mathcal{P}_3
\end{equation}
with $G_d$ and $\mathcal{P}_3$ as defined in {\rm (\ref{Gddef})} and
{\rm (\ref{poldef})}.
\end{theorem}

\vskip 0.8truecm

\begin{center}
\setlength{\unitlength}{0.35cm}
\begin{picture}(20,8)(2,0)
\put(0,0){\line(8,0){8}}
\put(0,0){\line(0,7){7}}
{\thicklines
\qbezier(0,5)(3,5)(6,5)
}
\put(-.8,-1.3){$0$}
\put(-1.5,4.8){$\gamma$}
\put(9,-.3){$\kappa$}
\put(-1,8){$\lambda_p(\kappa)$}
\put(15,0){\line(8,0){8}}
\put(15,0){\line(0,7){7}}
{\thicklines
\qbezier(15,4)(18,2.5)(21,2.2)
}
\qbezier[60](15,2)(18,2)(22,2)
\qbezier[60](15,5)(18,5)(22,5)
\put(14.2,-1.3){$0$}
\put(13.5,4.8){$\gamma$}
\put(13.5,1.8){$\rho\gamma$}
\put(24,-.3){$\kappa$}
\put(14,8){$\lambda_p(\kappa)$}
\put(0,5){\circle*{.45}}
\put(15,4){\circle*{.45}}
\end{picture}
\end{center}
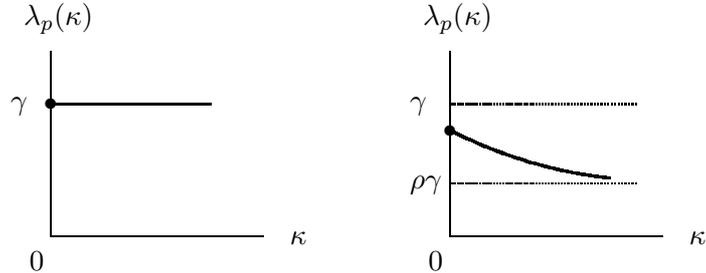
\caption{Qualitative picture of $\kappa\mapsto\lambda_p(\kappa)$
for recurrent, respectively, transient random walk.}\label{fig-SEqual}

\vskip 0.8truecm

\subsection{Discussion}
\label{S3.3}

The intuition behind Theorem~\ref{LyalowSE} is the following. If the catalyst is
driven by a recurrent random walk, then it suffers from ``traffic jams'', i.e.,
with not too small a probability there is a large region around the origin that
the catalyst fully occupies for a long time. Since with not too small a
probability the simple random walk (driving the reactant) can stay inside
this large region for the same amount of time, the average growth rate of the
reactant at the origin is maximal. This phenomenon may be expressed by saying
that \emph{for recurrent random walk clumping of the catalyst dominates the
growth of the moments}. For transient random walk, on the other hand, clumping
of the catalyst is present (the growth rate of the reactant is $>\rho\gamma$), 
but it is \emph{not\/} dominant (the growth rate of the reactant is $<\gamma$). 
Again, when the reactant stands still or moves slowly, the successive moments 
of the reactant are sensitive to successive degrees of clumping of the catalyst. 
As the diffusion constant $\kappa$ of the reactant increases, the effect of the 
clumping of the catalyst on the reactant gradually diminishes and the growth 
rate of the reactant gradually decreases to $\rho\gamma$.

Theorem~\ref{LyahighlimSE} has the same interpretation as its analogue
Theorem~\ref{LyapropIRW}(iii) for ISRW. We conjecture that the same behavior
occurs for SEP as in Conjectures~\ref{cj1IRW}--\ref{cj2IRW} for ISRW.


\section{Symmetric Voter Model}
\label{S4}

In this section we consider the case where $\xi$ is the \emph{Symmetric Voter Model}
(SVM) in equilibrium, or converging to equilibrium from a product measure. We summarize 
the results obtained in G\"artner, den Hollander and Maillard~\cite{garholmai07prb}.

\subsection{Model}
\label{S4.1}

As in Section~\ref{S3}, we abbreviate $\Omega=\{0,1\}^{\mathbb{Z}^d}$ and we let
$p\colon\mathbb{Z}^d\times\mathbb{Z}^d\to [0,1]$ be the transition kernel of an
irreducible symmetric random walk. The SVM is the Markov process on $\Omega$ whose
generator $L$ acts on cylindrical functions $f$ as
\begin{equation}
\label{VMgen}
(Lf)(\eta) = \sum_{x,y\in\mathbb{Z}^d} 1_{\{\eta(x)\neq\eta(y)\}}\,
p(x,y)\,\left[f(\eta^y)- f(\eta)\right],
\end{equation}
where
\begin{equation}
\label{VMtr}
\eta^y(z) =
\begin{cases}
\eta(z)   &\text{if } z \neq y,\\
1-\eta(y) &\text{if } z=y.
\end{cases}
\end{equation}
In words, site $x$ imposes its state on site $y$ at rate $p(x,y)$. The states
$0$ and $1$ are referred to as opinions or, alternatively, as vacancy and particle.
Contrary to ISRW and SEP, SVM is a non-conservative and non-reversible dynamics: 
opinions are not preserved.

We will consider two choices for the starting measure of $\xi$:
\begin{equation}
\label{muchoice}
\begin{cases}
&\nu_\rho, \mbox{ the Bernoulli product measure with density } \rho\in (0,1),\\
&\mu_\rho, \mbox{ the equilibrium measure with density } \rho\in (0,1).\\
\end{cases}
\end{equation}
The ergodic properties of the SVM are qualitatively different for recurrent and
for transient transition kernels. In particular, when $p(\cdot,\cdot)$ is recurrent
all equilibria are trivial, i.e., $\mu_\rho = (1-\rho)\delta_0+\rho\delta_1$,
while when $p(\cdot,\cdot)$ is transient there are also non-trivial equilibria,
i.e., ergodic $\mu_\rho$ parameterized by the density $\rho$. When starting from
$\nu_\rho$, $\xi(\cdot\,,t)$ converges in law to $\mu_\rho$ as $t\to\infty$.

\subsection{Main theorems}
\label{S4.2}

\begin{theorem}
\label{LyadiVM}
Let $d\geq 1$, $\kappa\in [0,\infty)$, $\rho\in (0,1)$, $\gamma\in (0,\infty)$ 
and $p\in\mathbb{N}$.\\
(i) For all $\kappa\in [0,\infty)$, the limit in {\rm (\ref{lyapdef})} exists
and is finite, and is the same for the two choices of starting measure in
$(\ref{muchoice})$.\\
(ii) On $\kappa\in [0,\infty)$, $\kappa\to\lambda_p(\kappa)$ is continuous.
\end{theorem}

The following dichotomy holds (see Fig.~\ref{fig-VMqual}):

\begin{theorem}
\label{LyalowVM}
Suppose that $p(\cdot,\cdot)$ has finite variance. Fix $\rho\in (0,1)$, $\gamma
\in (0,\infty)$ and $p\in\mathbb{N}$.\\
(i) If $1\leq d\leq 4$, then $\lambda_p(\kappa)=\gamma$ for all $\kappa \in 
[0,\infty)$.\\
(ii) If $d\geq 5$, then $\rho\gamma<\lambda_p(\kappa)<\gamma$ for all $\kappa
\in [0,\infty)$.
\end{theorem}

\vskip 0.8truecm

\begin{center}
\setlength{\unitlength}{0.35cm}
\begin{picture}(20,8)(0,0)
\put(0,0){\line(8,0){8}}
\put(0,0){\line(0,7){7}}
{\thicklines
\qbezier(0,5)(3,5)(6,5)
}
\put(-.8,-1.3){$0$}
\put(-1.5,4.8){$\gamma$}
\put(9,-.3){$\kappa$}
\put(-1,8){$\lambda_p(\kappa)$}
\put(15,0){\line(8,0){8}}
\put(15,0){\line(0,7){7}}
{\thicklines
\qbezier(15,4)(18,2.5)(21,2.2)
}
\qbezier[60](15,2)(18,2)(22,2)
\qbezier[60](15,5)(18,5)(22,5)
\put(14.2,-1.3){$0$}
\put(13.5,4.8){$\gamma$}
\put(13.5,1.8){$\rho\gamma$}
\put(24,-.3){$\kappa$}
\put(14,8){$\lambda_p(\kappa)$}
\put(0,5){\circle*{.45}}
\put(15,4){\circle*{.45}}
\end{picture}
\end{center}
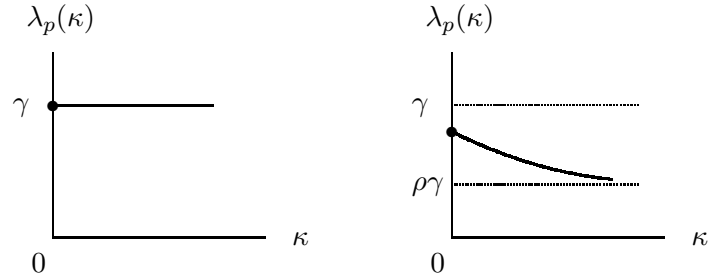
\caption{Qualitative picture of $\kappa\mapsto\lambda_p(\kappa)$ for symmetric 
random walk with finite variance in $d=1,2,3,4$, respectively, $d\geq 5$.}\label{fig-VMqual}

\vskip 0.8truecm

\begin{theorem}
\label{LyaintVM}
Suppose that $p(\cdot,\cdot)$ has finite variance. Fix $\rho\in (0,1)$ and $\gamma
\in (0,\infty)$. If $d\geq 5$, then $p\mapsto\lambda_p(0)$ is strictly increasing.
\end{theorem}

\subsection{Discussion}
\label{S4.3}

Theorem \ref{LyalowVM} shows that the Lyapunov exponents exhibit a dichotomy similar
to those found for ISRW and SEP (see Fig.~\ref{fig-VMqual}). The crossover in
dimensions is at $d=5$ rather than at $d=3$. Theorem \ref{LyaintVM} shows that
the system is intermittent at $\kappa=0$ when the Lyapunov exponents are nontrivial,
which is similar as well.

We conjecture that the following properties hold, whose analogues for ISRW and SEP
are known to be true:
\begin{conjecture}
\label{cj-monconvlim}
On $[0,\infty)$, $\kappa\mapsto\lambda_p(\kappa)$ is strictly decreasing and convex
with $\lim_{\kappa\to\infty}\lambda_p(\kappa)=\rho\gamma$.
\end{conjecture}

We close with a conjecture about the scaling behavior for $\kappa\to\infty$.
\begin{conjecture}
\label{cj-larg-kappaVM}
Let $d\geq 5$, $\rho \in (0,\infty)$ and $p\in\mathbb{N}$. Assume {\rm (\ref{SRW})}.
Then
\begin{equation}
\label{limlamb*VM}
\lim_{\kappa\to\infty} 2d\kappa[\lambda_p(\kappa)-\rho\gamma]
= \rho(1-\rho)\gamma^2 \frac{G_d^\ast}{G_d}
+ 1_{\{d=5\}}(2d)^3\left[\rho(1-\rho)\gamma^2\frac{1}{G_d}p\right]^2\mathcal{P}_5
\end{equation}
with
\begin{equation}
\label{Gdefs}
\begin{aligned}
G_d &= \int_0^\infty p_t(0,0)\, dt,\\
G_d^\ast &= \int_0^\infty t\, p_t(0,0)\, dt,
\end{aligned}
\end{equation}
and
\begin{equation}
\label{Rpdefshift}
\mathcal{P}_5 = \sup_{{f \in H^1(\mathbb{R}^5)} \atop {\|f\|_2=1}}
\Big[\, \int_{\mathbb{R}^5} dx\,|f(x)|^2 \int_{\mathbb{R}^5} dy\,|f(y)|^2
\,\frac{1}{16\pi^2\|x-y\|}
- \int_{\mathbb{R}^5} dx\,|\nabla f(x)|^2\,\Big].
\end{equation}
\end{conjecture}

\section{Concluding remarks}
\label{S5}

The theorems listed in Sections~\ref{S2}--\ref{S4} show that the intermittent
behavior of the reactant for the three types of catalyst exhibits interesting
similarities and differences. ISRW, SEP and SVM each show a dichotomy of strongly
catalytic versus weakly catalytic behavior, for ISRW between divergence and
convergence of the Lyapunov exponents, for SEP and SVM between maximality and
non-maximality. Each also shows an interesting dichotomy in the dimension for
the scaling behavior at large diffusion constants, with $d=3$ being critical 
for ISRW and SEP, and $d=5$ for SVM. For ISRW and SEP the same polaron term 
appears in the scaling limit, while for SVM an analogous but different polaron-like \
term appears. Although the techniques we use for the three models differ substantially, 
there is a universal principle behind their scaling behavior. See the heuristic 
explanation offered in \cite{garhol06} and \cite{garholmai07}. 

Both ISRW and SEP are conservative and reversible dynamics. The reversibility
allows for the use of spectral techniques, which play a key role in the analysis.
The SVM, on the other hand, is a non-conservative and irreversible dynamics.
The non-reversibility precludes the use of spectral techniques, and this dynamics
is therefore considerably harder to handle.

Both for SEP and SVM, the graphical representation is a powerful tool. For SEP
this graphical representation builds on random walks, for SVM on coalescing
random walks (see Liggett~\cite{lig85}).

The reader is invited to look at the original papers for details.



\end{document}